\documentclass[article,twoside,final,letterpaper,reqno,11pt]{amsart}
\usepackage{amssymb}
\usepackage{amsfonts}
\usepackage{amscd}

\newtheorem{theorem}{Theorem}
\numberwithin{theorem}{subsubsection}
\makeatletter
\def \fixequation{\let \c@equation\c@theorem\let \p@equation\p@theorem\let \theequation \thetheorem}
\def \@doit#1{\numberwithin{theorem}{#1}\fixequation}
\let \@tempstartsection \@startsection
\def \@startsection#1#2#3#4#5#6{\ifx#1\@empty \else \@doit{#1}\fi \@tempstartsection{#1}{#2}{#3}{#4}{#5}{#6}}
\makeatother
\theoremstyle{plain}

\newtheorem{definition}[theorem]{Definition}

\newtheorem{lemma}[theorem]{Lemma}

\newtheorem{proposition}[theorem]{Proposition}
\newtheorem{remark}[theorem]{Remark}

\makeatletter
\let \p@equation\p@theorem
\let \c@equation\c@theorem 
\let \theequation \thetheorem 
\makeatother
\input{tcilatex}
\begin{document}
\title{The Geometric Weil Representation}
\author{Shamgar Gurevich}
\address{Department of Mathematics, University of California, Berkeley, CA
94720-3840, USA.}
\email{shamgar@math.berkeley.edu}
\author{Ronny Hadani}
\address{Department of Mathematics, University of Chicago, Chicago, IL
60637, USA.}
\email{hadani@math.uchicago.edu}
\date{February 1, 2006.}
\thanks{\copyright \ Copyright by S. Gurevich and R. Hadani, February 2006.
All rights reserved.}

\begin{abstract}
In this paper we construct a geometric analogue of the Weil representation
over a finite field. Our construction is principally invariant, not choosing
any specific realization. This eliminates most of the unpleasant formulas
that appear in the traditional (non-invariant) approaches, and puts in the
forefront some delicate geometric phenomena which underlie this
representation.
\end{abstract}

\maketitle

\section{Introduction}

\subsection{The Weil representation}

In his celebrated 1964 Acta paper \cite{W} Andr\`{e} Weil constructed a
certain (projective) unitary representation of a symplectic group over local
fields. This representation has many fascinating properties which have
gradually been brought to light over the last few decades. It now appears
that this representation is a central object in mathematics, bridging
between various topics in mathematics and physics, including classical
invariant theory, the theory of theta functions and automorphic forms,
harmonic analysis, and last (but probably not least) quantum mechanics.
Although it holds such a fundamental status, it is satisfying to observe
that the Weil representation already appears in "real life" situations.
Given a linear space $L$, there exists an associated (split) symplectic
vector space $V=L\times L^{\ast }$. The Weil representation of the group $%
Sp=Sp(V,\omega )$ can be realized on the Hilbert space 
\begin{equation*}
\mathcal{H=}\text{ }L^{2}(L,%
\mathbb{C}
).
\end{equation*}

Interestingly, elements of the group $Sp$ act by certain kinds of
generalized Fourier transforms. In particular, there exists a specific
element\ $\mathrm{w}\in Sp$ \ (called the Weyl element) whose action is
given, up to a normalization, by the standard Fourier transform. From this
perspective, the classical theory of harmonic analysis seems to be devoted
to the study of a particular (and not very special) operator in the Weil
representation.

In this paper we will be concerned only with the case of the Weil
representations of symplectic groups over finite fields of odd
characteristic. Let us note that already in this setting there is no simple
way to obtain the Weil representation. A possible approach to its
construction is through the use of an auxiliary two-step unipotent group
called the \textit{Heisenberg group}. In this approach, the Weil
representation appears as a collection of intertwining operators of a
special irreducible representation of the Heisenberg\textit{\ }group. It
takes some additional work to realize that in fact, in the finite field
case, it can be linearized into an honest representation.

\subsection{Geometric approach and invariant presentation}

The Weil representation over a finite field is an object of algebraic
geometry! More specifically, in the case of finite fields, all groups
involved are finite, yet of a very special kind, consisting of rational
points of corresponding algebraic groups. In this setting it is an ideology,
due to Grothendieck, that any (meaningful) set-theoretic object is governed
by a more fundamental algebro-geometric one. The procedure by which one
lifts to the setting of algebraic geometry is called \textit{geometrization. 
}In this procedure, sets are replaced by algebraic varieties (defined over
the finite field) and functions are replaced by corresponding
sheaf-theoretic objects (i.e., objects of the $\ell $-adic derived
category). The procedure translating between the two settings is called 
\textit{sheaf-to-functions correspondence. }This ideology has already proved
to be extremely powerful, not just as a technique for proving theorems, but
much more importantly, it supplies an appropriate framework without which\
certain deep mathematical ideas cannot even be stated.

It is reasonable to suspect that the Weil representation is meaningful
enough so that it can be geometrized. Indeed, this was carried out in a
letter written by Deligne to Kazhdan in 1982 \cite{D1}. In his letter
Deligne proposed a geometric analogue of the Weil representation. More
precisely, he proposed to construct an $\ell $-adic Weil perverse sheaf on
the algebraic group $\mathbf{Sp}$ (in this paper we use boldface letters to
denote algebraic varieties).

Although the main ideas already appear in that letter, the content is
slightly problematic for several reasons. First, it was never published; and
this is probably a good enough reason for writing this text. Second, in his
letter, Deligne chooses a specific realization \cite{Ge}, which is
equivalent to choosing a splitting $V=L\times L^{\prime }$. As a result he
obtains a sheaf%
\begin{equation*}
\mathcal{K}_{\mathbf{L,L}^{\prime }}\text{ on the variety }\mathbf{Sp\times L%
}^{\prime }\mathbf{\times L}^{\prime }\mathbf{.}
\end{equation*}%
The sheaf \ $\mathcal{K}_{\mathbf{L,L}^{\prime }}$ is given in terms of
explicit formulas which of course depend on the choice of $L$ and $L^{\prime
}$. These arbitrary choices make the formulas quite unpleasant (see \cite{GH}
for the $SL_{2}$ case) and, moreover, help to hide some delicate geometric
phenomena underlying the Weil representation.

In this paper, we present a construction of the geometric Weil
representation in all dimensions, yet from a different perspective. In our
construction we do not use any specific realization. Instead, we invoke the
idea of \textit{invariant presentation} of an operator and as a consequence
we obtain an invariant construction of the Weil representation sheaf%
\begin{equation*}
\mathcal{K}\text{ \ on the variety }\mathbf{Sp\times V.}
\end{equation*}%
This avoids most of the unpleasant computations, and moreover, brings to the
forefront some of the geometry behind this representation. As a by-product
of working in the geometric setting, we obtain a new proof of the linearity
of the Weil representation over finite fields. The upshot is that it is
enough to prove the multiplicativity property on an open subvariety of the
symplectic group. It is possible to choose such a variety on which verifying
multiplicativity becomes simple. This kind of argument has no analogue in
the set-theoretic setting and it is our hope that it exemplifies some of the
powerful elegance of the geometric method.

Although it was first proposed some 25 years ago, we would like to remark
that recently a striking application was found for the geometric Weil
representation in the area of quantum chaos and the associated theory of
exponential sums \cite{GH}. We believe that this is just the tip of the
iceberg and the future promises great things for this fundamental object.
This belief is our motivation to write this paper.

\subsection{Results}

\begin{enumerate}
\item The main result of this paper is the new invariant construction of the
geometric Weil representation sheaf in all dimensions. This is the content
of Theorem \ref{GOR} which claims the existence and some of the properties
of this sheaf.

\item We present a new construction of the set-theoretic Weil representation
in all dimensions, which uses the notion of invariant presentation of an
operator. In addition, algebro-geometric techniques (such as perverse
extension) are used to obtain a new proof of the multiplicativity (i.e.,
linearity) property. Finally, explicit formulas are supplied.

\item We present a character formula for the Weil representation. This is
the content of Theorem \ref{CF}.
\end{enumerate}

\subsection{Structure of the paper}

The paper is divided into several sections.

\begin{itemize}
\item Section \ref{preliminaries}. In this section we present some
preliminaries. In the first part of this section, the Heisenberg
representation is presented and the notion of invariant presentation of an
operator is discussed. In the second part of this section the classical
construction of the Weil representation is given.

\item Section \ref{IPW}. The Weil representation and its character are
obtained using the idea of invariant presentation of an operator.

\item Section \ref{geometricoscillator}. The geometric Weil representation
sheaf is constructed. The construction is obtained by means of
Grothendieck's geometrization procedure concluding with Theorem \ref{GOR},
which claims the existence of a sheaf-theoretic analogue of the Weil
representation. Some preliminaries from algebraic geometry and $\ell $-adic
sheaves are also discussed. In particular, we recall Grothendieck's
sheaf-to-function\/ correspondence procedure which is the main connection
between Sections \ref{geometricoscillator} and \ref{IPW}.

\item Appendix \ref{DS}. The relation between our invariant construction of
the Weil representation sheaf and the sheaf that was proposed by Deligne 
\cite{D1} is discussed. The upshot is that they are related by a certain
kind of $\ell $-adic Fourier transform.

\item Appendix \ref{GT}. A proof of the main theorem of the paper, i.e.,
Theorem \ref{GOR} is given.

\item Appendix \ref{P}. Certain additional proofs are given
\end{itemize}

\subsubsection{ Acknowledgments}

It is a pleasure to thank our teacher J. Bernstein for his guidance and for
informing us of Deligne's letter. We thank D. Kazhdan for several important
remarks. We thank A. Weinstein for the discussion on the Weyl transform. We
appreciate the efforts of O. Gabber and S. Lysenko who read the first draft
of this paper. We would like to acknowledge P. Deligne for allowing us to
publish his ideas about the geometrization of the Weil representation. In
addition, we are grateful to R. Howe for sharing with us unpublished results
of his concerning the Cayley transform. We thank V. Drinfeld for the
opportunity to present this paper in the "Geometric Langlands Seminar". We
appreciate the support of M. Baruch. We would like to thank J. Keating,
Bristol, and the IHES, Paris, for the hospitality during our visit in
February 2006, where part of this work was done.

\section{Preliminaries\label{preliminaries}}

In this paper we denote by $k=\mathbb{F}_{q}$ a finite field of $q$
elements, where $q$ is odd, and by $\overline{k}$ an algebraic closure.

\subsection{The Heisenberg representation}

Let $(V,\omega )$ be a $2N$-dimensional symplectic vector space over the
finite field $k$. There exists a two-step nilpotent group $H=H\left(
V,\omega \right) $ associated to the symplectic vector space $(V,\omega )$.
\ The group $H$ is called the \textit{Heisenberg group.} It can be realized
as the set $H=V\times k$ equipped with the following multiplication rule 
\begin{equation}
(v,z)\cdot (v^{\prime },z^{\prime })=(v+v^{\prime },z+z^{\prime }+\tfrac{1}{2%
}\omega (v,v^{\prime })).  \label{mult}
\end{equation}

The center of $H$ is $Z=Z(H)=\{(0,z):z\in k\}$. $\ $Fix a non-trivial
additive character $\psi :Z\longrightarrow 
\mathbb{C}
^{\times }$. We have the following fundamental theorem

\begin{theorem}[Stone-Von Neumann]
\label{SVN}

There exists a unique $($up to isomorphism$)$ irreducible representation $%
(\pi ,H,\mathcal{H})$ with central character $\psi $, i.e., $\pi (z)=\psi
(z)Id_{\mathcal{H}}$ for every $z\in Z.$
\end{theorem}

We call the representation $\pi ,$ appearing in Theorem \ref{SVN}, the 
\textit{Heisenberg representation} associated with the central character $%
\psi $.

\begin{remark}
\label{Sch}Although the representation $\pi $ is unique, it admits a
multitude of different models $($realizations$)$. In fact this is one of its
most interesting and powerful attributes. In particular, for any Lagrangian
splitting $V=L\oplus L^{\prime }$, there exists the Schr\"{o}dinger model $%
(\pi _{L,L^{\prime }},H,S(L^{\prime }))$, where $S(L^{\prime })$ denotes the
space of complex-valued functions on $L^{\prime }$. In this model we have
the following actions

\begin{itemize}
\item $[\pi _{L,L^{\prime }}(l)\vartriangleright f](x)=\psi (\omega
(x,l))f(x);$

\item $[\pi _{L,L^{\prime }}(l^{\prime })\vartriangleright
f](x)=f(x+l^{\prime })$;

\item $[\pi _{L,L^{\prime }}(z)\vartriangleright f](x)=\psi (z)f(x)$%
,\smallskip
\end{itemize}
\end{remark}

where $l\in L,$ $x,$ $l^{\prime }\in L^{\prime }$, and $z\in Z$.

\subsection{Invariant presentation of operators\label{IP}}

Let $A:\mathcal{H\longrightarrow H}$ be a linear operator. Choosing an
identification $\beta $ of $\mathcal{H}$ with the space $S(X)$ of functions
on a finite set $X$ (which is equivalent to choosing a basis of $\mathcal{H)}
$, we can realize $A$ as a function $A_{\beta }\in S(X\times X)$ (which is
equivalent to writing $A$ as a matrix). The function $A_{\beta }$ acts by
convolution on functions in $S(X)$. Choosing a different identification $%
\beta ,$ of course, changes the realization $A_{\beta }$. This is all very
trivial, nevertheless, we would like to pose the following question:

\begin{center}
\textbf{Is there a correct realization?}
\end{center}

Representation theory suggests an answer to this seemingly strange question.
We will realize operators as certain functions on the group $H$. In more
detail, associated to a central character $\psi :Z\longrightarrow 
\mathbb{C}
^{\times }$ we denote by $S(H,\psi )$ the space of functions on the
Heisenberg group which behave $\psi ^{-1}$-equivariantly with respect to
left (right) multiplication of the center, i.e., functions $f\in S(H)$
satisfying 
\begin{equation}
f(z\cdot h)=\psi ^{-1}(z)f(h),  \label{EP}
\end{equation}%
for every $z\in Z$.

Recall that $\mathcal{H}$ constitutes an irreducible representation space of
the Heisenberg group (with central character $\psi $). This implies that
given an operator $A$ there exists a unique function $\widehat{A}\in
S(H,\psi )$ such that $\pi (\widehat{A})=A$, where $\pi (\widehat{A})=\frac{1%
}{|Z(H)|}\tsum \limits_{h\in H}\widehat{A}(h)\pi (h)$. The relation between $%
A $ and $\widehat{A}$ is expressed in the following formula 
\begin{equation}
\widehat{A}(h)=\tfrac{1}{dim\mathcal{H}}Tr(A\pi (h^{-1})).
\label{presentation}
\end{equation}%
The transform $A\longmapsto \widehat{A}$ is well known in the literature and
is usually referred to as the \textit{Weyl transform \cite{We2}}. It behaves
well with respect to composition of operators: Given two operators $A,B$ $%
\in End(\mathcal{H})$, it is a direct calculation to verify that $\widehat{AB%
}=\widehat{A}\ast \widehat{B}$, where $\ast $ denotes the normalized
convolution product of functions on the group

\begin{equation}
\widehat{A}\ast \widehat{B}(h)=\tfrac{1}{|Z(H)|}\tsum \limits_{h_{1}\cdot
h_{2}=h}\widehat{A}(h_{1})\widehat{B}(h_{2}).  \label{convolution}
\end{equation}

\bigskip

Concluding, using the Weyl transform we are able to establish an invariant
presentation of operators.

\begin{theorem}[Invariant presentation of operators \protect \cite{H2}]
Let $(\pi ,H,\mathcal{H})$ be the Heisenberg representation. The Weyl
transform $\widehat{\left( \cdot \right) }:End(\mathcal{H})\rightarrow
S(H,\psi )$ is an isomorphism of algebras, with the inverse given by the
extended action $\pi :$ $S(H,\psi )\rightarrow End(\mathcal{H})$, 
\begin{equation*}
f\mapsto \pi (f)=\frac{1}{|Z(H)|}\tsum \limits_{h\in H}f(h)\pi (h).
\end{equation*}
\end{theorem}

Because of the equivariance property (\ref{EP}), functions in $S(H,\psi )$
are determined by their restriction to the vector space $V\subset H,$ hence
we can (and will) identify $S(H,\psi )$ with $S(V)$. The Convolution (\ref%
{convolution}) is realized on $S(V)$ as follows, for two functions $f,$ $g$ $%
\in S(V)$

\begin{equation*}
(f\ast g)(v)=\tsum \limits_{v_{1}+v_{2}=v}\psi (-\tfrac{1}{2}\omega
(v_{1},v_{2}))f(v_{1})g(v_{2}).
\end{equation*}%
We will refer to the induced isomorphism $\widehat{\left( \cdot \right) }%
:End(\mathcal{H})\rightarrow S(V)$ also as the Weyl transform.

\subsection{The Weil representation \label{WR}}

Let $G=Sp(V,\omega )$ be the group of symplectic linear automorphisms of $V$%
. The group $G$ acts by group automorphism on the Heisenberg group through
its tautological action on the vector space $V$. This induces an action of $%
G $ on the category $Rep(H)$ of representations of $H$, i.e., given a
representation $\pi \in Rep(H)$ and an element $g\in G$ one obtains a new
representation $\pi ^{g}$ (realized on the same Hilbert space) defined by $%
\pi ^{g}(h)=\pi (g\cdot h)$. It is clear that this action does not affect
the central character and sends an irreducible representation to an
irreducible one. Let $\pi $ be the Heisenberg representation associated with
a central character $\psi $. Invoking Theorem \ref{SVN}, we conclude that
for every element $g\in G$ we have%
\begin{equation*}
\pi ^{g}\simeq \pi .
\end{equation*}%
Denote by $\rho (g):\mathcal{H\longrightarrow H}$ an intertwiner which
realizes this isomorphism$.$ Equivalently, this means that $\rho (g)$
satisfies and in fact is determined up to scalar by the following equation

\begin{equation}
\rho (g)\pi (h)\rho (g)^{-1}=\pi (g\cdot h),  \label{Egorov2}
\end{equation}%
for every $g\in G.$

The above equation are sometimes referred to in the literature as the Egorov
identity. Having that all $\pi ^{g},$ $g\in G,$ are irreducibles, and using
Schur's lemma, we conclude that the collection $\{ \rho (g);$ $g\in G\}$
forms a projective representation $\rho :G\longrightarrow PGL(\mathcal{H}).$
It is a non-trivial argument that the projective representation $\rho $ can
be linearized into an honest representation which we also denote by $\rho
:G\longrightarrow GL(\mathcal{H})$. This representation is called the 
\textit{Weil representation}. Let us summarize this in the following theorem

\begin{theorem}
\label{linearization}There exists a canonical\footnote{%
Unique unless $q=3$ and $N=1$ $($see Remark \ref{q3N1} for the canonical
choice in the latter case$).$} representation 
\begin{equation*}
\rho :G\longrightarrow GL(\mathcal{H)},
\end{equation*}%
satisfying the equation (\ref{Egorov2})$.$
\end{theorem}

\begin{remark}
The fact that the projective representation $\rho $ can be linearized is a
peculiar phenomenon of the finite field situation. If one deals with the
analogue constructions, for example over infinite Archimedian fields, it is
a deep fact that $\rho $ can be de-projectivized up to a $\pm $ sign, which
is called the metaplectic sign. Hence, in this case $\rho $ can be corrected
to a representation of a double cover%
\begin{equation*}
1\rightarrow 
\mathbb{Z}
_{2}\rightarrow \widetilde{Sp}\rightarrow Sp\rightarrow 1,
\end{equation*}%
which is called the metaplectic cover. This fact is responsible for several
fundamental phenomena in quantum mechanics, such as the nonzero energy of
the vacuum state of the harmonic oscillator.
\end{remark}

\section{Invariant presentation of the Weil representation \label{IPW}}

The Weil representation was obtained classically by a very implicit
construction (Section \ref{WR}). This is already true before the
linearization procedure (Theorem \ref{linearization}). What we would like to
do next is to obtain a concrete realization (formulas) for the
representation $\rho $. In order to carry this out, it seems that one would
be obliged to choose a basis of the Hilbert space $\mathcal{H}$ and then
write every element $\rho (g)$ as an appropriate matrix. This procedure, in
principle, could be carried out and, in fact, in most references we know
this is the way $\rho $ is presented. We would like to suggest a different
approach based on the idea of \textit{invariant presentation} discussed in
Section\  \ref{IP}. Namely, every operator $\rho (g)$ will be presented by a
function (kernel) $K_{g}=\widehat{\rho (g)}\in S(V)$ satisfying the
multiplicative rule $K_{gh}=K_{g}\ast K_{h}$. The collection of kernels $%
\{K_{g}\}_{g\in G}$ \ can be thought of as a single function%
\begin{equation}
K\in S(G\times V),  \label{WKF}
\end{equation}%
satisfying the following convolution property 
\begin{equation}
(m\times Id)^{\ast }K=p_{1}^{\ast }K\ast p_{2}^{\ast }K,
\label{convolution2}
\end{equation}%
where the map $m\times Id:G\times G\times V\longrightarrow G\times V$ is
given by $(g_{1},g_{2,}v)\mapsto (g_{1}\cdot g_{2},v)$, and the maps $%
p_{i}(g_{1},g_{2,}v)=(g_{i},v)$ are the projectors on the first and second $%
G $-coordinate respectively.\ The right-hand side of (\ref{convolution2})
means \ $p_{1}^{\ast }K\ast p_{2}^{\ast }K(g_{1},g_{2},v)=K_{g_{1}}\ast
K_{g_{2}}(v)$. In the sequel we will sometime suppress the $V$-coordinate,
writing $m:G\times G\longrightarrow G$, and $p_{i}:G\times G\longrightarrow
G $, we will also suppress the projections $p_{i}$ in (\ref{convolution2})
obtaining a much cleaner convolution formula%
\begin{equation}
m^{\ast }K=K\ast K.  \label{cc}
\end{equation}

\subsection{Ansatz\label{An}}

We would like to establish an ansatz for the kernels $K_{g}$ for elements $g$
in some appropriate subset of $G$. Let $O$ be the subset of all elements $%
g\in End(V)$ such that \bigskip $g-I$ is invertible. We define the following
map $\kappa :O\longrightarrow End(V)$%
\begin{equation*}
\kappa (g)=\frac{g+I}{g-I}.
\end{equation*}%
It is a direct verification that $\kappa $ is injective, identifying the set 
$U=G\cap O$ with a subset of the Lie algebra $sp(V)$. The map $\kappa $ is
well known in the literature \cite{We1} and referred to as the \textit{%
Cayley transform}. It establishes an interesting algebraic relation between
the group $G$ and its Lie algebra.

We summarize some basic properties of the Cayley transform

\begin{itemize}
\item $\kappa (g^{-1})=-\kappa (g).$ \ 

\item $\kappa ^{2}=Id.$

\item $\omega (\kappa (g)v,u)=-\omega (v,\kappa (g)u)$ for every $g\in U$
(this is just the statement that $\kappa (g)\in sp(V)).$

\item The following two equivalent identities hold%
\begin{eqnarray}
\kappa (gh) &=&[I+\kappa (g)][\kappa (g)+\kappa (h)]^{-1}[I-\kappa
(g)]+\kappa \left( g\right) ,  \label{iden1_eq} \\
\kappa (gh) &=&[I+\kappa (h)][\kappa (g)+\kappa (h)]^{-1}[I+\kappa (g)]-I,
\label{iden2_eq}
\end{eqnarray}

for every pair $g,h\in U$ such that $gh\in U$. We note that for such a pair $%
\kappa (g)+\kappa (h)$ is invertible, therefore, the formulas indeed make
sense.
\end{itemize}

Let $\sigma :k^{\times }\longrightarrow 
\mathbb{C}
^{\times }$ be the unique quadratic character of the multiplicative group
(also called the Legendre character). Denote by $\mathfrak{e}$ the value of
the one-dimensional Gauss sum $\mathfrak{e=}\tsum \limits_{z\in Z}\psi
(z^{2}) $.

\begin{description}
\item[Ansatz] Define%
\begin{equation}
K_{g}(v)=\nu (g)\psi (\tfrac{1}{4}\omega (\kappa (g)v,v)),\text{ \  \  \ }g\in
U,  \label{anzats}
\end{equation}%
where $\nu (g)=\tfrac{\mathfrak{e}^{2N}}{q^{2N}}\sigma (\det (\kappa
(g)+I)). $
\end{description}

\begin{remark}
The only non-trivial part of the ansatz (\ref{anzats}) is the normalization
factor $\nu (g)$. The expression which involves the $\psi $ can be easily
obtained \cite{H2}. In more detail, the operator $\rho (g)$ is determined up
to a scalar by the Egorov identity (\ref{Egorov2}). Hence, starting with the
identity operator $Id:\mathcal{H\longrightarrow H}$ and taking the average \ 
\begin{equation*}
A=\tsum \limits_{v\in V}\pi (-g\cdot v)Id\pi (v)=\pi (\widehat{A})\in
Hom_{H}(\pi ,\pi ^{g}\mathcal{)},
\end{equation*}%
we recover the operator $\rho (g)$ up to a scalar. It is a direct
calculation, using formula $($\ref{presentation}$)$, to verify that $%
\widehat{A}(v)=\psi (\tfrac{1}{4}\omega (\kappa (g)v,v))$.
\end{remark}

\begin{theorem}[Extension]
\label{ET}There exists a unique multiplicative extension of \ (\ref{anzats})
to $G.$
\end{theorem}

The uniqueness is clear since we have $U\cdot U=G$, the existence of this
extension will be proved in Section \ref{geometricoscillator} using
algebraic geometry.

\begin{remark}
In the case where $q=3$ and $N=1$ the Weil representation is not unique, the
extension established in Theorem \ref{ET} gives, in a sense that will be
explained later (see Remark \ref{q3N1}), a natural choice$.$
\end{remark}

\subsection{Application to characters}

Denote by $\tau =\rho \ltimes \pi $ the representation of the semi-direct
product $J=G\ltimes H$ on the Hilbert space $\mathcal{H}$. The group $J$ is
sometime referred to as the Jacobi group. We will call the representation $%
(\tau ,J,\mathcal{H)}$ the \textit{Heisenberg-Weil representation}. As a
result of the explicit expression (\ref{anzats}) we obtain, using formula (%
\ref{presentation}), a new (cf. \cite{Ge,H1, H2}) explicit formulas for the
character of the Weil and Heisenberg-Weil representations when restricted to
the subsets $U$ and $U\times H,$ respectively.

\begin{theorem}[Character formulas]
\label{CF}The character $ch_{\tau \text{ }}$of the Heisenberg-Weil
representation $\tau $ when restricted to the subset $U\times H$ is given by 
\begin{equation*}
ch_{\tau }(g,v,z)=\tfrac{\mathfrak{e}^{2N}}{q^{N}}\sigma (\det (\kappa
(g)+I))\psi (\tfrac{1}{4}\omega (\kappa (g)v,v)+z),
\end{equation*}%
and the character $ch_{\rho \text{ }}$of the Weil representation $\rho $
when restricted to the subset $U$ is given by 
\begin{equation*}
ch_{\rho }(g)=\tfrac{\mathfrak{e}^{2N}}{q^{N}}\sigma (\det (\kappa (g)+I)).
\end{equation*}
\end{theorem}

\section{Geometric Weil representation\label{geometricoscillator}}

In this section we are going to construct a geometric counterpart to the
set-theoretic Weil representation. This will be achieved using \textit{%
geometrization}. This is a formal procedure, invented by Grothendieck and
his school, by which sets are replaced by algebraic varieties (over the
finite field) and functions are replaced by certain sheaf-theoretic objects.

\subsection{Preliminaries from algebraic geometry}

Next we have to use some space to recall notions and notations from
algebraic geometry and the theory of $\ell $-adic sheaves. \ 

\subsubsection{Varieties}

In the sequel, we are going to translate back and forth between algebraic
varieties defined over the finite field $k$ and their corresponding sets of
rational points. In order to prevent confusion between the two, we use
bold-face letters to denote a variety $\mathbf{X}$ and normal letters $X$ to
denote its corresponding set of rational points $X=\mathbf{X}(k)$. For us, a
variety $\mathbf{X}$ over the finite field is a quasi-projective algebraic
variety, such that the defining equations are given by homogeneous
polynomials with coefficients in the finite field $k$. In this situation,
there exists a (geometric) \textit{Frobenius} endomorphism $Fr:\mathbf{%
X\rightarrow X}$, which is a morphism of algebraic varieties. We denote by $%
X $ \ the set of points fixed by $Fr$, i.e., $X=\mathbf{X}(k)=\mathbf{X}%
^{Fr}=\{x\in \mathbf{X}:Fr(x)=x\}.$ The category of algebraic varieties over 
$k$ will be denoted by $\mathsf{Var}_{k}$.

\subsubsection{Sheaves}

Let $\mathsf{D}^{b}(\mathbf{X)}$ denote the bounded derived category of
constructible $\ell $-adic sheaves on $\mathbf{X}$ \cite{BBD}. We denote by $%
\mathsf{Perv}(\mathbf{X)}$ the Abelian category of perverse sheaves on the
variety $\mathbf{X}$, that is the heart with respect to the autodual
perverse t-structure in $\mathsf{D}^{b}(\mathbf{X})$. An object $\mathcal{%
F\in }\mathsf{D}^{b}(\mathbf{X)}$ is called $n$-perverse if $\mathcal{F[}%
n]\in \mathsf{Perv}(\mathbf{X)}$. Finally, we recall the notion of a Weil
structure (Frobenius structure) \cite{D2}. A Weil structure associated to an
object $\mathcal{F\in }\mathsf{D}^{b}(\mathbf{X)}$ is an isomorphism%
\begin{equation*}
\theta :Fr^{\ast }\mathcal{F}\overset{\sim }{\longrightarrow }\mathcal{F}%
\text{.}
\end{equation*}

A pair $(\mathcal{F},\theta )$ is called a Weil object. By an abuse of
notation we often denote $\theta $ also by $Fr$. We choose once an
identification $\overline{%
\mathbb{Q}
}_{\ell }\simeq 
\mathbb{C}
$, hence all sheaves are considered over the complex numbers.

\begin{remark}
All the results in this section make perfect sense over the field $\overline{%
\mathbb{Q}
}_{\ell }$, in this respect the identification of $\overline{%
\mathbb{Q}
}_{\ell }$ with $%
\mathbb{C}
$ \ is redundant. The reason it is specified is in order to relate our
results with the standard constructions.
\end{remark}

Given a Weil object $(\mathcal{F},Fr^{\ast }\mathcal{F\simeq F})$ one can
associate to it a function $f^{\mathcal{F}}:X\rightarrow 
\mathbb{C}
$ to $\mathcal{F}$ as follows 
\begin{equation*}
f^{\mathcal{F}}(x)=\tsum \limits_{i}(-1)^{i}Tr(Fr_{|H^{i}(\mathcal{F}_{x})}).
\end{equation*}%
This procedure is called \textit{Grothendieck's sheaf-to-function
correspondence}. Another common notation for the function $f^{\mathcal{F}}$
is $\chi _{Fr}(\mathcal{F)}$, which is called the \textit{Euler
characteristic} of the sheaf $\mathcal{F}.$

\subsection{Geometric Weil representation}

We shall now start the geometrization of the Weil representation.

\subsubsection{Replacing sets by varieties}

The first step we take is to replace all sets involved by their geometric
counterparts, i.e., algebraic varieties. The symplectic space $(V,\omega )$
is naturally identified as the set $V=\mathbf{V}(k)$, where $\mathbf{V}$ is
a $2N$-dimensional symplectic vector space in $\mathsf{Var}_{k}$. The
Heisenberg group $H$ is naturally identified as the set $H=\mathbf{H}(k)$,
where $\mathbf{H}=\mathbf{V\times }\mathbb{A}^{1}$\ is the group variety
equipped with the same multiplication formulas (\ref{mult}). The subvariety $%
\mathbf{Z}=Z(\mathbf{E})=\{(0,\lambda ):\lambda \in \mathbb{A}^{1}\}$ is the
center of $\mathbf{H}$. Finally, the group $G$ is naturally identified as $G=%
\mathbf{G}(k)$, where $\mathbf{G}=Sp(\mathbf{V,}\omega )$.

\subsubsection{Replacing functions by sheaves}

The second step is to replace functions by their sheaf-theoretic
counterparts \cite{Ga}. The central character $\psi :Z\longrightarrow 
\mathbb{C}
^{\times }$ is associated via the sheaf-to-function correspondence to the
Artin-Schreier sheaf $\mathcal{L}_{\psi }$ living on $\mathbf{Z}$, i.e., we
have $f^{\mathcal{L}_{\psi }}=\psi .$ The Legendre character $\sigma $ on $%
k^{\times }\simeq $ $\mathbb{G}_{m}(k)$ is associated to the Kummer sheaf $%
\mathcal{L}_{\sigma }$ on $\mathbb{G}_{m}$. Looking back at formula (\ref%
{anzats}), the constant $\mathfrak{e}$ is replaced by the following Weil
object 
\begin{equation*}
\mathcal{E=}\tint \limits_{\mathbf{Z}}\mathcal{L}_{\psi (z^{2})}\in \mathsf{D%
}^{b}(\mathbf{pt}),
\end{equation*}%
where, for the rest of this paper, $\int =\int_{!}$ denotes integration with
compact support \cite{BBD}. Grothendieck's Lefschetz trace formula \cite{Gr}
implies that, indeed, $f^{\mathcal{E}}=\mathfrak{e}.$ In fact, there exists
a quasi-isomorphism $\mathcal{E}\overset{q.i}{\longrightarrow }$ $H^{1}(%
\mathcal{E)}[-1]$ and $\dim H^{1}(\mathcal{E)=}1$, hence, $\mathcal{E}$ can
be thought of as a one-dimensional vector space, equipped with a Frobenius
operator, sitting at cohomological degree $1.$

Finally, given an element $a$ in the Lie algebra $sp(\mathbf{V)}$, we will
use the notation%
\begin{equation*}
\mathcal{G}_{a}=\tint \limits_{v\in \mathbf{V}}\mathcal{L}_{\psi (\frac{1}{4}%
\omega (av,v))},
\end{equation*}%
for the symplectic gauss integral (noting that $\omega (a\cdot ,\cdot )$ is
a symmetric bilinear form).

Our main objective, in this section, is to construct a multiplicative
extension $K:G\times V\longrightarrow 
\mathbb{C}
$ of the anzats (\ref{anzats}). The extension appears as a direct
consequence of the following geometrization theorem (see Appendix \ref{GT}
for a proof).

\begin{theorem}[Geometric Weil representation]
\label{GOR} There exists a geometrically irreducible $\dim (\mathbf{G)}$%
-perverse Weil sheaf $\mathcal{K}$ on $\mathbf{G\times V}$ of pure weight $w(%
\mathcal{K})=0$, satisfying the following properties

\begin{enumerate}
\item Convolution property:\label{property1} There exist an isomorphism 
\begin{equation}
m^{\ast }\mathcal{K\simeq K\ast K}.  \label{CP}
\end{equation}

\item Function property:\label{property2} We have $f^{\mathcal{K}}\left(
g,v\right) =K_{g}\left( v\right) $ for every $g\in U.$
\end{enumerate}
\end{theorem}

\bigskip

\begin{remark}
Note that, as promised in Theorem \ref{ET}, taking $K=f^{\mathcal{K}}$, we
obtain a multiplicative extension of the ansatz! The nice thing about this
construction is that it uses geometry and, in particular (see Appendix \ref%
{CO}), the notion of perverse extension which has no counterpart in the
set-function theoretic setting.
\end{remark}

\subsection{Additional properties}

We describe two additional properties of the sheaf $\mathcal{K}$.

\subsubsection{Restriction property}

We consider an irreducible subvariety $\mathbf{M\subset G}$ and the
associated open subvariety $\mathbf{U}_{\mathbf{M}}=\mathbf{U\cap M}$, where 
$\mathbf{U}$ $\mathbf{=\{}g\mathbf{\in G;}$ $g-I$ is invertible$\mathbf{\}.}$

\begin{definition}
We say that $\mathbf{M}$ is \underline{multiplicative} if $\mathbf{M\cdot
M\subset M.}$ A multiplicative subvariety $\mathbf{M\subset G}$ is called 
\underline{openly generated} if the induced morphism%
\begin{equation}
m:\mathbf{U}_{\mathbf{M}}\mathbf{\times U}_{\mathbf{M}}\rightarrow \mathbf{M,%
}  \label{OG}
\end{equation}%
is smooth, surjective and with connected fibers.
\end{definition}

\begin{theorem}
\label{RP} Let $\mathbf{M\subset G}$ be a multiplicative, openly generated
subvariety, then the restriction $\mathcal{K}_{|\mathbf{M\times V}}$ is
geometrically irreducible $\dim (\mathbf{M)}$-perverse Weil sheaf of pure
weight zero.
\end{theorem}

For a proof see Appendix \ref{PRP}.

\subsubsection{Product property}

Consider two symplectic vector spaces $\mathbf{\left( V_{1},\omega
_{1}\right) }$ and $\mathbf{\left( V_{2},\omega _{2}\right) }$ in $\mathsf{%
Var}_{k}$. Denote by $\mathbf{\left( V,\omega \right) }$ the product space
with $\mathbf{V=V}_{1}\times \mathbf{V}_{2}$ and $\mathbf{\omega =\omega }%
_{1}+\mathbf{\omega }_{2}$. Let $\mathbf{G}_{1,}\mathbf{G}_{2}$ and $\mathbf{%
G}$ be the corresponding symplectic groups. We have the natural inclusion $%
i: $ $\mathbf{G}_{1}\times \mathbf{V}_{1}\times \mathbf{G}_{2}\times \mathbf{%
V}_{2}\hookrightarrow \mathbf{G\times V}$. Let us denote by $\mathcal{K}_{1},%
\mathcal{K}_{2}$ and $\mathcal{K}$ the corresponding Weil representation
sheaves.

\begin{theorem}
\label{Pp}There exist an isomorphism $i^{\ast }\mathcal{K\simeq K}%
_{1}\boxtimes $ $\mathcal{K}_{2}$, where $\boxtimes $ denotes exterior
tensor product.
\end{theorem}

For a proof see Appendix \ref{PPp}.

\begin{remark}[\protect \cite{G}]
\label{q3N1}The product property confirms that the set theoretic Weil
representation that we constructed coincides with the one in \cite{Ge}. This
is the unique family $\left( \rho _{V},\text{ }Sp\left( V\right) ,\mathcal{H}%
_{V}\right) $, where\ $V$ runs in the category of symplectic vector spaces,
which satisfies 
\begin{equation*}
\rho _{V_{1}\times V_{2}}\simeq \rho _{V_{1}}\boxtimes \rho _{V_{2}}\text{,}
\end{equation*}%
for every pair of symplectic vector spaces $V_{1}$ and $V_{2}$.
\end{remark}

\appendix

\section{Relation to Deligne's sheaf\label{DS}}

In his letter \cite{D1} Deligne proposed a sheaf-theoretic analogue of the
Weil representation of $G=Sp(V,\omega )$, where $\left( V,\omega \right) $
is a $2N$-dimensional symplectic vector space over the finite field $k=%
\mathbb{F}_{q}$ with $q$ odd$.$ His proposed construction was given in terms
of specific formulas. In this short section we would like to elaborate on
the relation between the Weil representation sheaf $\mathcal{K}$ which we
constructed in the previous subsection, and Deligne's sheaf (see \cite{GH}
for the formal construction of Deligne's sheaf in the $SL_{2}$ case).

Deligne's sheaf was given in terms of a specific realization of the Weil
representation. Let us elaborate on his construction. Fixing a Lagrangian
splitting $V=L\oplus L^{\prime },$ we consider the Schr\"{o}dinger
realization (Remark \ref{Sch}) $(\pi _{L,L^{\prime }},H,\mathcal{H}%
_{L,L^{\prime }})$ of the Heisenberg representation, where $\mathcal{H}%
_{L,L^{\prime }}=S(L^{\prime })$. In this realization every operator $\rho
(g)$, $\ g\in G,$ is given by a kernel $K_{L,L^{\prime }}(g):L^{\prime
}\times L^{\prime }\rightarrow 
\mathbb{C}
$. The collection of kernels $\left \{ K_{L,L^{\prime }}(g)\right \} _{g\in
G}$ is equivalent \ to a single function $K_{L,L^{\prime }}$ on $G\times
L^{\prime }\times L^{\prime }$ which satisfies the convolution property%
\begin{equation}
m^{\ast }K_{L,L^{\prime }}=K_{L,L^{\prime }}\ast K_{L,L^{\prime }}.
\label{convolution4}
\end{equation}

In his letter Deligne constructed a perverse Weil sheaf, which we will
denote by $\mathcal{K}_{\mathbf{L,L}^{\prime }}$, on the variety $\mathbf{%
G\times L}^{\prime }\times \mathbf{L}^{\prime }$ satisfying \ the following
two properties

\begin{itemize}
\item Function property\textbf{: } The function $K_{L,L^{\prime }}=$ $f^{%
\mathcal{K}_{\mathbf{L,L}^{\prime }}}$.

\item Convolution property:\textbf{\ } There exists an isomorphism%
\begin{equation*}
m^{\ast }\mathcal{K}_{\mathbf{L,L}^{\prime }}\simeq \mathcal{K}_{\mathbf{L,L}%
^{\prime }}\ast \mathcal{K}_{\mathbf{L,L}^{\prime }}.
\end{equation*}
\end{itemize}

Let us now explain the relation between our sheaf $\mathcal{K}$ and
Deligne's sheaf $\mathcal{K}_{\mathbf{L,L}^{\prime }}$. On the level of
functions, $K$ and $K_{L,L^{\prime }}$ are related as follows%
\begin{eqnarray*}
K_{L,L^{\prime }}(g,x,y) &=&\left \langle \delta _{x}|\rho (g)\delta
_{y}\right \rangle \\
&=&\tsum \limits_{l^{\prime }\in L^{\prime },l\in L}K(g,l^{\prime }+l)\left
\langle \delta _{x}|\psi (\omega (\cdot ,l)L_{l^{\prime }}\delta _{y}\right
\rangle \psi (\tfrac{1}{2}\omega (l,l^{\prime })) \\
&=&\tsum \limits_{l^{\prime }\in L^{\prime },l\in L}K(g,l^{\prime }+l)\left
\langle \delta _{x}|\delta _{y-l^{\prime }}\right \rangle \psi (\omega
(y-l^{\prime },l)-\tfrac{1}{2}\omega (l,l^{\prime })) \\
&=&\tsum \limits_{l\in L}K(g,y-x+l)\psi (\omega (x,l)-\tfrac{1}{2}\omega
(l,y-x)),
\end{eqnarray*}%
where the first equality is the definition of $K_{L,L^{\prime }}(g,x,y)$ and
the second equality follows from substituting $\pi _{L,L^{\prime
}}(l^{\prime },l)=\psi (\tfrac{1}{2}\omega (l,l^{\prime }))M_{\psi (\omega
(\cdot ,l)}L_{l^{\prime }}$. Concluding, we obtained%
\begin{equation}
K_{L,L^{\prime }}(g,x,y)=\tsum \limits_{l\in L}K(g,y-x+l)\psi (\tfrac{1}{2}%
\omega (x+y,l)).  \label{formula13}
\end{equation}%
It will be convenient to rewrite (\ref{formula13}) in a more functorial
manner. Let $\alpha :G\times L^{\prime }\times L^{\prime }\rightarrow
G\times L^{\prime }\times L^{\prime }$ be the map given by $\alpha
(g,x,y)=(g,y-x,x+y).$ Formula (\ref{formula13}) is equivalent to%
\begin{equation}
K_{L,L^{\prime }}=\alpha ^{\ast }Four_{L}(K),  \label{formula14}
\end{equation}%
where $Four_{L}$ denotes the (non-normalized) fiber-wise Fourier transform
from functions on the vector bundle $G\times L^{\prime }\times L\rightarrow
G\times L^{\prime }$ to functions on the dual vector bundle $G\times
L^{\prime }\times L^{\prime }\rightarrow G\times L^{\prime }$; the duality
is manifested by the symplectic form $\omega $ which induces a
non-degenerate pairing between $L$ and $L^{\prime }$.

Formula (\ref{formula14}) can easily be translated into the geometric
setting. We define Deligne's sheaf associated with the Schr\"{o}dinger
realization $\pi _{L,L^{\prime }}$ to be%
\begin{equation*}
\mathcal{K}_{\mathbf{L,L}^{\prime }}=\alpha ^{\ast }Four_{L}(\mathcal{K}),
\end{equation*}%
where now $Four_{L}$ denotes the $\ell $-adic Fourier transform \cite{KL}.
Since $Four_{L}$ shifts the perversity degree by $N$ \cite{KL} it follows
that $\mathcal{K}_{\mathbf{L,L}^{\prime }}$ is an irreducible $\mathrm{\dim }%
(\mathbf{G)}+N$-perverse Weil sheaf.

\section{Proof of the geometric Weil representation theorem \label{GT}}

Appendix \ref{GT} is devoted to the proof of Theorem \ref{GOR}, i.e., to the
construction of the Weil representation sheaf $\mathcal{K}$.

\subsection{Construction\label{CO}}

The construction of the sheaf $\mathcal{K}$ is based on the ansatz (\ref%
{anzats}).

Let $\mathbf{U\subset G}$ be the open subvariety consisting of elements $%
g\in \mathbf{G}$ such that $g-I$ is invertible. The Cayley transform is
given by an algebraic function, hence it make sense also in the geometric
setting, we denote it also by $\kappa .$ The construction proceeds in
several steps

\begin{itemize}
\item \textit{Non-normalized kernels:}\textbf{\  \ }On the variety $End\left( 
\mathbf{V}\right) \mathbf{\times V}$ define the sheaf%
\begin{equation*}
\widetilde{\mathcal{K}}\left( a,v\right) =\mathcal{L}_{\psi (\tfrac{1}{4}%
\omega (av,v))}.
\end{equation*}

Denote by $\widetilde{\mathcal{K}}_{\mathbf{U}}$ the pull back $\widetilde{%
\mathcal{K}}_{\mathbf{U}}\left( g,v\right) =$ $\widetilde{\mathcal{K}}\left(
\kappa \left( g\right) ,v\right) $ defined on $\mathbf{U\times V}$.

\item \textit{Normalization coefficients:}\textbf{\ }On the open subvariety $%
\mathbf{U\times V}$ define the sheaf%
\begin{equation*}
\mathcal{C=E}^{\otimes 2N}\otimes \mathcal{L}_{\sigma \left( \det (\kappa
(g)+I)\right) }[4N]\left( 2N\right) .
\end{equation*}

We note that for $g\in \mathbf{U}$ we have $\kappa (g)+I=\frac{2g}{g-I}$,
thus it is invertible and consequently $\mathcal{L}_{\sigma \left( \det
(\kappa (g)+I)\right) }$ is smooth on $\mathbf{U}$.

\item \textit{Normalized kernels:}\textbf{\ }On the open subvariety $\mathbf{%
U\times V}$ define the sheaf%
\begin{equation}
\mathcal{K}_{\mathbf{U}}=\mathcal{C\otimes }\widetilde{\mathcal{K}}_{\mathbf{%
U}}.  \label{NK}
\end{equation}%
Finally, take%
\begin{equation}
\mathcal{K=}j_{!\ast }\mathcal{K}_{\mathbf{U}},  \label{kernelsheaf}
\end{equation}%
where $j:\mathbf{U\hookrightarrow G}$ is the open imbedding, and $j_{!\ast 
\text{ }}$is the functor of perverse extension \cite{BBD}. It follows
directly from the construction that the sheaf $\  \mathcal{K}$ is irreducible 
$\dim (\mathbf{G)}$-perverse of pure weight $0$.
\end{itemize}

Section \ref{PP} is devoted to proving that the kernel sheaf $\mathcal{K\ }$(%
\ref{kernelsheaf}) satisfies the conditions of Theorem \ref{GOR}.

\subsection{Proof of the properties\label{PP}}

The function property is clear from the construction. We are left to show
the convolution property.

We need to show that%
\begin{equation}
m^{\ast }\mathcal{K\simeq K\ast K},  \label{property1'}
\end{equation}%
where $m:\mathbf{G\times G\longrightarrow G}$ is the multiplication morphism.

\begin{proof}
The morphism $m$ forms a principal $\mathbf{G}$-bundle, hence it is smooth
and surjective. We denote by $m_{\mathbf{U}}:\mathbf{U\times
U\longrightarrow G}$ the restriction of $m$ to the open subvariety $\mathbf{%
U\times U\hookrightarrow G\times G}$. It is easy to verify that $m_{\mathbf{U%
}}$ is smooth and surjective. The following technical lemma (see Appendix %
\ref{PTL} for a proof) plays a central role in the proof.

\begin{lemma}
\label{technicallemma} There exists an isomorphism $m_{\mathbf{U}}^{\ast }%
\mathcal{K\simeq K}_{\mathbf{U}}\mathcal{\ast K}_{\mathbf{U}}$. \smallskip
\smallskip
\end{lemma}

Lemma \ref{technicallemma} implies that the sheaves $m^{\ast }\mathcal{K}$
and $\mathcal{K\ast K}$ are isomorphic when restricted to the open
subvariety $\mathbf{U\times U}$. The sheaf $m^{\ast }\mathcal{K}$ is
irreducible $\dim (\mathbf{G\times G)}$-perverse as a pull-back by a smooth
morphism of an irreducible $\dim (\mathbf{G)}$-perverse sheaf on $\mathbf{%
G\times V}$. Hence, it is enough to show that the sheaf $\mathcal{K\ast K}$
is also irreducible $\dim (\mathbf{G\times G)}$-perverse. Consider the map $%
(m_{\mathbf{U}},Id):\mathbf{U\times U\times G\longrightarrow G\times G}$. It
is clearly smooth and surjective. It is enough to show that the pull-back $%
(m_{\mathbf{U}},Id)^{\ast }\mathcal{K\ast K}$ is irreducible $\dim (\mathbf{%
U\times U\times G)}$-perverse. Again, using Lemma \ref{technicallemma} and
also invoking some direct diagram chasing one obtains%
\begin{equation}
(m,Id)^{\ast }\mathcal{K\ast K\simeq K}_{\mathbf{U}}\ast \mathcal{K}_{%
\mathbf{U}}\ast \mathcal{K}.  \label{formula1}
\end{equation}%
The right-hand side of (\ref{formula1}) is principally a subsequent
application of a properly normalized, symplectic Fourier transforms (see
Formula (\ref{Fo}) below) on $\mathcal{K}$, hence by the Katz-Laumon theorem 
\cite{KL} it is irreducible $\dim (\mathbf{U\times U\times G)}$-perverse.

Let us summarize. We showed that both sheaves $m^{\ast }\mathcal{K}$ and $%
\mathcal{K\ast K}$ are irreducible $\dim (\mathbf{G\times G)}$-perverse and
are isomorphic on an open subvariety. This implies that they must be
isomorphic. This concludes the proof of Property $($\ref{property1'}$)$.
\end{proof}

\subsection{ Proof of Lemma \protect \ref{technicallemma}\label{PTL}}

The proof will proceed in several steps.

\begin{itemize}
\item \bigskip Step 1. We show that the sheaf $\mathcal{K}_{\mathbf{U}}\ast 
\mathcal{K}_{\mathbf{U}}$ is irreducible $\dim (\mathbf{U\times U)}$%
-perverse. This is done by a direct computation%
\begin{equation*}
\mathcal{K}_{\mathbf{U}}\ast \mathcal{K}_{\mathbf{U}}(g,h,v)\simeq \mathcal{C%
}_{g}\otimes \mathcal{C}_{h}\otimes \widetilde{\mathcal{K}}_{\mathbf{U}}\ast 
\widetilde{\mathcal{K}}_{\mathbf{U}}(g,h,v).
\end{equation*}%
Writing the convolution more explicitly\ 
\begin{equation}
\widetilde{\mathcal{K}}_{\mathbf{U}}\ast \widetilde{\mathcal{K}}_{\mathbf{U}%
}(g,h,v)\simeq \widetilde{\mathcal{K}}_{\mathbf{U}}(g,v)\otimes \tint
\limits_{u\in \mathbf{V}}\mathcal{L}_{\psi (\tfrac{1}{2}\omega (v,[\kappa
(g)-I]u)}\otimes \widetilde{\mathcal{K}}(\kappa (g)+\kappa (h),u).
\label{Fo}
\end{equation}

Hence, we see that $\widetilde{\mathcal{K}}_{\mathbf{U}}\ast \widetilde{%
\mathcal{K}}_{\mathbf{U}}$ is principally an application of Fourier
transform, which implies that $\mathcal{K}_{\mathbf{U}}\ast \mathcal{K}_{%
\mathbf{U}}$ is irreducible $\dim (\mathbf{U\times U)}$-perverse.

\item Step 2. It is enough to show that the sheaves $m^{\ast }\mathcal{K}$
and $\mathcal{K}_{\mathbf{U}}\ast \mathcal{K}_{\mathbf{U}}$ are isomorphic
on an open subvariety. Let $\mathbf{W\subset U\times U}$ be the open
subvariety consisting of pairs $(g,h)\in \mathbf{U\times U}$ such that $%
gh\in \mathbf{U}$. We will show that the sheaves $m^{\ast }\mathcal{K}$ and $%
\mathcal{K}_{\mathbf{U}}\ast \mathcal{K}_{\mathbf{U}}$ are isomorphic on $%
\mathbf{W}$. Direct computation of the integral in (\ref{Fo}), using
completion to squares, reveals that%
\begin{equation*}
\  \widetilde{\mathcal{K}}_{\mathbf{U}}\ast \widetilde{\mathcal{K}}_{\mathbf{U%
}}(g,h,v)\simeq \widetilde{\mathcal{K}}(\kappa (g)+[I+\kappa (g)][\kappa
(g)+\kappa (h)]^{-1}[I-\kappa (g)],v)\otimes \mathcal{G}_{\kappa (g)+\kappa
(h)},
\end{equation*}%
for every $\left( g,h\right) \in W$.
\end{itemize}

Invoking identity (\ref{iden1_eq}), we obtain that%
\begin{equation*}
\widetilde{\mathcal{K}}_{\mathbf{U}}\ast \widetilde{\mathcal{K}}_{\mathbf{U}%
}(g,h,v)\simeq \widetilde{\mathcal{K}}_{\mathbf{U}}(gh,v)\otimes \mathcal{G}%
_{\kappa (g)+\kappa (h)}.
\end{equation*}

Consequently, it is enough to show

\begin{proposition}
\label{isom}There exists an isomorphism 
\begin{equation*}
\mathcal{C}_{g}\otimes \mathcal{C}_{h}\otimes \mathcal{G}_{\kappa (g)+\kappa
(h)}\simeq \mathcal{C}_{gh}.
\end{equation*}
\end{proposition}

\begin{proof}
Denote by $T$ the Tate vector space $T=H^{2}(\tint \limits_{\mathbb{P}^{1}}%
\mathbb{C}
_{\mathbb{P}^{1}})[-2].$ Recall that we defined $\mathcal{C}_{g}=\mathcal{E}%
^{\otimes 2N}\otimes \mathcal{L}_{\sigma (\det [\kappa (g)+I])}\left[ 4N%
\right] \left( 2N\right) $ and $\mathcal{C}_{h}=\mathcal{E}^{\otimes
2N}\otimes \mathcal{L}_{\sigma (\det [\kappa (h)+I])}\left[ 4N\right] \left(
2N\right) .$ Using the character property of the sheaf $\mathcal{L}_{\sigma }
$ we obtain%
\begin{eqnarray}
\mathcal{C}_{g}\otimes \mathcal{C}_{h} &\simeq &\mathcal{E}^{\otimes
4N}\otimes \mathcal{L}_{\sigma (\det [\kappa (g)+I]\det [\kappa (h)+I])}%
\left[ 8N\right] \left( 4N\right)   \label{formula4} \\
&\simeq &\{ \mathcal{E}^{\otimes 4}\}^{N}\otimes \mathcal{L}_{\sigma (\det
[\kappa (g)+I]\det [\kappa (h)+I])}\left[ 8N\right] \left( 4N\right)   \notag
\\
&\simeq &\{T^{\otimes 2}\}^{\otimes N}\otimes \mathcal{L}_{\sigma (\det
[\kappa (g)+I]\det [\kappa (h)+I])}\left[ 8N\right] \left( 4N\right)   \notag
\\
&\simeq &T^{\otimes 2N}\otimes \mathcal{L}_{\sigma (\det [\kappa (g)+I]\det
[\kappa (h)+I]).}\left[ 8N\right] \left( 4N\right)   \notag \\
&\simeq &\mathcal{L}_{\sigma (\det [\kappa (g)+I]\det [\kappa (h)+I])}\left[
4N\right] (2N).
\end{eqnarray}%
In addition, short analysis of the symplectic Gauss integral reveals%
\begin{equation}
\mathcal{G}_{\kappa (g)+\kappa (h)}\simeq \mathcal{E}^{\otimes 2N}\otimes 
\mathcal{L}_{\sigma (\det [\kappa (g)+\kappa (h)]).}  \label{formula5}
\end{equation}%
Combining formulas (\ref{formula4}), (\ref{formula5}), and identity (\ref%
{iden2_eq}), we obtain%
\begin{equation*}
\mathcal{C}_{g}\otimes \mathcal{C}_{h}\otimes \mathcal{G}_{\kappa (g)+\kappa
(h)}\simeq \mathcal{E}^{\otimes 2N}\otimes \mathcal{L}_{\sigma (\det [\kappa
(gh)+I])}[4N]\left( 2N\right) =\mathcal{C}_{gh}.
\end{equation*}
\end{proof}

\section{Additional proofs\label{P}}

\subsection{Proof of Theorem \protect \ref{RP}\label{PRP}}

\begin{proof}
Denote $\mathcal{F=K}_{|\mathbf{M\times V}}$ and $\mathcal{F}_{\mathbf{U}}=%
\mathcal{K}_{|\mathbf{U}_{M}\mathbf{\times V}}$.

Since restriction does not increase weight, it is enough to prove that $%
\mathcal{F}$ is geometrically irreducible $\dim (\mathbf{M)}$-perverse. Due
to the assumption that $m$ is smooth and surjective, it is enough to show
that $m^{\ast }\mathcal{F}$ (where $m$ is the morphism (\ref{OG})) is a
geometrically irreducible $\dim (\mathbf{U}_{\mathbf{M}}\times \mathbf{U}_{%
\mathbf{M}})$-perverse. \smallskip We have $m^{\ast }\mathcal{F}\simeq 
\mathcal{F}_{\mathbf{U}}\ast \mathcal{F}_{\mathbf{U}}$, implying that $%
m^{\ast }\mathcal{F}$ is an application of a properly normalized $\ell $%
-adic Fourier transform to an irreducible smooth sheaf (See formula (\ref{Fo}%
)), consequently, the assertion follows by Katz-Laumon theorem \cite{KL}%
.\smallskip
\end{proof}

\subsection{Proof of Theorem \protect \ref{Pp}\label{PPp}}

\begin{proof}
\textit{Step 1. }The sheaves $i^{\ast }\mathcal{K}$\ and $\mathcal{K}%
_{1}\boxtimes \mathcal{K}_{2}$ are geometrically irreducible $[\dim (\mathbf{%
G}_{1}\times \mathbf{G}_{2})]$-perverse Weil sheaves on $\mathbf{G}%
_{1}\times \mathbf{G}_{2}$. The statement is evident for $\mathcal{K}%
_{1}\boxtimes \mathcal{K}_{2},$ and for $i^{\ast }\mathcal{K}$ it follows
from the restriction property (Theorem \ref{RP}) of the sheaf $\mathcal{K}$
with respect to the multiplicative subvariety $\mathbf{M=G}_{1}\times 
\mathbf{G}_{2}\subset $ $\mathbf{G}$.

\textit{Step 2.}\textbf{\ }By perversity it is enough to show that $i^{\ast }%
\mathcal{K}$\ and $\mathcal{K}_{1}\boxtimes \mathcal{K}_{2}$ are isomorphic
when restricted to the open subvariety $\mathbf{U}_{1}\times \mathbf{V}%
_{1}\times \mathbf{U}_{2}\times \mathbf{V}_{2},$ where $\mathbf{U}%
_{i}=\{g\in \mathbf{G}_{i};$ $\det (g-I)\neq 0\}$. This is a simple
verification using formula (\ref{NK}).
\end{proof}

\end{document}